\documentclass{article}
\pdfoutput=1
\usepackage{amsmath}
\usepackage{amsfonts}
\usepackage{epsfig}
\usepackage{rotating}
\usepackage{graphicx}
\usepackage{amssymb}
\usepackage{eucal}
\usepackage{lscape}
\usepackage{latexsym,enumerate}
\usepackage{caption}
\usepackage{setspace}
\usepackage{fancyhdr}
\usepackage{pdfpages}
\usepackage{natbib}
\usepackage{moreverb}
\usepackage{xr-hyper}

\usepackage[colorlinks,citecolor=blue,linkcolor=blue,urlcolor=blue]{hyperref}

\begin{document}
\title {Demystified: double robustness with nuisance parameters estimated at rate n-to-the-1/4}
\author{Judith J. Lok\\ Department of Mathematics and Statistics, Boston University\\ jjlok@bu.edu}

\maketitle

\section*{Abstract}

Have you also been wondering what is this thing with double robustness and nuisance parameters estimated at rate $n^{1/4}$? It turns out that to understand this phenomenon one just needs the Middle Value Theorem (or a Taylor expansion) and some smoothness conditions. This note explains why under some fairly simple conditions, as long as the nuisance parameter $\theta\in\mathbb{R}^k$ is estimated at rate $n^{1/4}$ or faster, 1.\ the resulting variance of the estimator of the parameter of interest $\psi\in\mathbb{R}^d$ does not depend on how the nuisance parameter $\theta$ is estimated, and 2.\ the sandwich estimator of the variance of $\hat{\psi}$ ignoring estimation of $\theta$ is consistent.

\section{Introduction}

It is not uncommon that an estimator for a parameter $\psi$ depends on nuisance parameters $\theta$. In such settings, $\theta$ is often estimated in a first step. Some estimators for $\psi$ are doubly robust: they depend on two nuisance parameters $\theta_1$ and $\theta_2$, and are consistent if one of the nuisance parameters $\theta_1$ or $\theta_2$ is consistently estimated, but not necessarily both.

Double robustness has been shown to often improve precision, and several efficient estimators that depend on more than one nuisance parameter have been shown to be doubly robust. Examples of this include doubly robust estimation of means from observational data (\cite{Bang}), doubly robust estimation of (coarse) Structural Nested Mean Models (\cite{comm,optLok}), and multiply robust estimation of indirect and direct effects (\cite{tchetgen2012semiparametric}). The orthogonal moment functions from \cite{chernozhukov2022locally} are locally doubly robust (see their equation~(2.4)), but beyond the scope of this note.

In order to obtain the efficiency gain from double robustness, it is advantageous to use flexible models to estimate $\theta$. Flexible methods do not always estimate $\theta$ at rate $\sqrt{n}$ (e.g., \cite{chernozhukov2022locally}). Fortunately, it often suffices to estimate $\theta$ at rate $n^{1/4}$ in order to obtain the efficiency gain, and if this is achieved, the variance of the resulting estimator $\hat{\psi}$ does not depend on how the nuisance parameter $\theta$ is estimated. 

It does not take much more than the Mean Value Theorem (or a Taylor expansion) to understand this phenomenon. This note shows how this works for estimators $\hat{\psi}$ based on smooth unbiased estimating equations. 

\section{Setting and Notation}\label{psisection}

Henceforth, $\psi^*\in\mathbb{R}^d$ is the true parameter of interest and $\theta^*\in\mathbb{R}^k$ is the true nuisance parameter. $\hat{\psi}$ solves
\begin{equation}\label{psihat}
\mathbb{P}_n U(\psi,\hat{\theta})=0,
\end{equation}
where $\mathbb{P}_n$ denotes the empirical average over $i=1,\ldots,n$ independent identically distributed observations, with
\begin{equation}\label{ident}
EU\bigl(\psi^*,\theta^*\bigr)=0
\end{equation}
and $U$ of dimension $k$, the dimension of $\psi$. Examples include Maximum Likelihood Estimation settings where $\hat{\psi}$ solves the score equations, but this so-called Z-estimation is much more general; see e.g.\ \cite{Vaart}. 

Such $\hat{\psi}$ is doubly robust if with $\theta=\bigl(\theta_1,\theta_2\bigr)$, $\hat{\psi}$ solves unbiased estimating equations if $\theta_1$ is consistently estimated and if $\theta_2$ is consistently estimated, and not necessarily both; that is,
\begin{equation}\label{DR}
EU(\psi^*,\theta^*_1,\theta_2)=0 \;\;\;\;\;\text{ and }\;\;\;\;\; EU(\psi^*,\theta_1,\theta^*_2)=0
\end{equation}
for every $\theta_1$ and $\theta_2$.

This note assumes that $\theta$ is estimated at rate $n^{1/4}$ or faster:
\begin{equation}\label{noneforth}
n^{1/4}\bigl(\hat{\theta}-\theta^*\bigr)=O_P(1).
\end{equation}

\section{Regularity conditions}\label{regconditions}

Throughout, this note assumes that regularity conditions hold so that $U\bigl(\psi,\theta\bigr)$ and $EU\bigl(\psi,\theta\bigr)$ depend smoothly enough on $\bigl(\psi,\theta\bigr)$. It also assumes that the order of differentiation and integration with respect to $\theta$ can be changed, so it is assumed that the support of the distribution of the observations does not depend on $\theta$. 

It is also assumed that $\psi$ is uniquely identified by equation~(\ref{ident}), so that
\begin{equation*}
\left.\frac{\partial}{\partial\psi}\right|_{\psi^*}EU\bigl(\psi,\theta^*\bigr)
\end{equation*}
has an inverse. To simplify the exposition, it is assumed that it has already been proven that $\hat{\psi}$ converges in probability to $\psi^*$.

\section{Derivations based on Taylor expansions}

Double robustness implies that
\begin{equation}\label{DRderivative}
E\left.\frac{\partial}{\partial\theta_p}\right|_{\theta^*} U_q(\psi^*,\theta)=0,
\end{equation}
where $U_q$ is the $q$th component of $U$, $1\leq q\leq k\in\mathbb{N}$.
This follows for the derivative with respect to $\theta_1$ by taking the derivative with respect to $\theta_1$ of $EU\bigl(\psi^*,\theta_1,\theta_2^*\bigr)$, which equals zero because of equation~(\ref{DR}). Notice that this assumes that the support of the observations does not depend on $\theta$, so that differentiation with respect to $\theta$ and integration can be interchanged. The same reasoning works for the derivative with respect to $\theta_2$.

After estimating $\theta$ resulting in $\hat{\theta}$, $\hat{\psi}$ solves equation~(\ref{psihat}):
\begin{eqnarray}0&=&\mathbb{P}_n U(\hat{\psi},\hat{\theta})\nonumber\\
&=&\mathbb{P}_n U\bigl(\psi^*,\theta^*\bigr)
+
\left(\left.\frac{\partial}{\partial\bigl(\psi,\theta\bigr)}\right|_{(\tilde{\psi},\tilde{\theta})}\mathbb{P}_n U\bigl(\psi,\theta\bigr)\right)
\left(\begin{array}{c}\hat{\psi}-\psi^*\\\hat{\theta}-\theta^*\end{array}\right)\nonumber\\
&=&\mathbb{P}_n U\bigl(\psi^*,\theta^*\bigr)
+
\left(\left.\frac{\partial}{\partial\psi}\right|_{(\tilde{\psi},\tilde{\theta})}\mathbb{P}_n U\bigl(\psi,\theta\bigr)\right)
\bigl(\hat{\psi}-\psi^*\bigr)\nonumber\\
&&+\left(\left.\frac{\partial}{\partial\theta}\right|_{(\tilde{\psi},\tilde{\theta})}\mathbb{P}_n U\bigl(\psi,\theta\bigr)\right)\bigl(\hat{\theta}-\theta^*\bigr)\nonumber\\
&=&\mathbb{P}_n U\bigl(\psi^*,\theta^*\bigr)
+
\left(\left.\frac{\partial}{\partial\psi}\right|_{\tilde{\psi}}\mathbb{P}_n U(\psi,\tilde{\theta})\right)
\bigl(\hat{\psi}-\psi^*\bigr)\nonumber\\
&&+\left(\left.\frac{\partial}{\partial\theta}\right|_{\tilde{\theta}}\mathbb{P}_n U\bigl(\tilde{\psi},\theta\bigr)\right)\bigl(\hat{\theta}-\theta^*\bigr)\nonumber\\
&&\label{starter}
\end{eqnarray}
for some $(\tilde{\psi},\tilde{\theta})$ between $(\hat{\psi},\hat{\theta})$ and $(\psi^*,\theta^*)$, possibly different in each row (from the Middle Value Theorem applied to each entry in the vector separately). Equation~(\ref{starter}) implies that
\begin{eqnarray}
\lefteqn{\bigl(\hat{\psi}-\psi^*\bigr)=}\nonumber\\
&&\left(\left.\frac{\partial}{\partial\psi}\right|_{\tilde{\psi}}\mathbb{P}_n U\bigl(\psi,\tilde{\theta}\bigr)\right)^{-1} \left(\mathbb{P}_n U\bigl(\psi^*,\theta^*\bigr)-\left.\frac{\partial}{\partial\theta}\right|_{\tilde{\theta}}\mathbb{P}_n U\bigl(\tilde{\psi},\theta\bigr)\bigl(\hat{\theta}-\theta^*\bigr)\right).\label{psiminpsistar}
\end{eqnarray}
The derivations below show that if equation~(\ref{noneforth}) holds,
the last term in equation~(\ref{psiminpsistar}) multiplied by $\sqrt{n}$ converges in probability to zero.

First, we show that equation~(\ref{noneforth}) implies that
\begin{equation}\label{psinforth}
n^{1/4}\bigl(\hat{\psi}-\psi^*\bigr)
\rightarrow^P0.
\end{equation}
Notice that as usual (see for example \cite{Lok} Lemma~A.6.1), under the usual regularity conditions (mainly differentiability conditions), since $\bigl(\tilde{\psi},\tilde{\theta}\bigr)\rightarrow^P \bigl(\psi^*,\theta^*\bigr)$,
\begin{equation}\label{der}\left.\frac{\partial}{\partial\psi}\right|_{\tilde{\psi}}\mathbb{P}_n U(\psi,\tilde{\theta})\rightarrow^{P} E\left.\frac{\partial}{\partial\psi}\right|_{\psi^*}U(\psi,\theta^*)
\end{equation}
and
\begin{equation}\left.\frac{\partial}{\partial\theta}\right|_{\tilde{\theta}}\mathbb{P}_n U\bigl(\tilde{\psi},\theta\bigr)\rightarrow^P E\left.\frac{\partial}{\partial\theta}\right|_{\theta^*} U(\psi^*,\theta)=0,
\end{equation}
where the equality follows from the double robustness equation~(\ref{DRderivative}).
Combining with equation~(\ref{noneforth}), it follows that the last term in equation~(\ref{psiminpsistar}) multiplied by $n^{1/4}$ converges in probability to zero. Combining with the Central Limit Theorem on $\mathbb{P}_n U\bigl(\psi^*,\theta^*\bigr)$, equation~(\ref{psiminpsistar}) shows that equation~(\ref{noneforth}) implies equation~(\ref{psinforth}).

To show that the last term in equation~(\ref{psiminpsistar}) multiplied by $\sqrt{n}$ converges in probability to zero, we next consider each
\begin{equation}\label{ddthetap}
n^{1/4}\left.\frac{\partial}{\partial\theta_p}\right|_{\tilde{\theta}}\mathbb{P}_n U_q\bigl(\tilde{\psi},\theta\bigr)
\end{equation}
separately, were $U_q$ is the $q$th component of $U$. We show that the quantity in equation~(\ref{ddthetap}) converges in probability to zero when equation~(\ref{noneforth}) holds.
\begin{equation}\label{doubleder}
\left.\frac{\partial}{\partial\theta_p}\right|_{\tilde{\theta}}\mathbb{P}_n U_q\bigl(\tilde{\psi},\theta\bigr)
=\left.\frac{\partial}{\partial\theta_p}\right|_{\theta^*}\mathbb{P}_n U_q\bigl(\tilde{\psi},\theta\bigr)+\left.\frac{\partial}{\partial\theta}\right|_{\dot{\theta}}\frac{\partial}{\partial\theta_p}\mathbb{P}_n U_q\bigl(\tilde{\psi},\theta\bigr)\bigl(\tilde{\theta}-\theta^*\bigr)
\end{equation}
because of the Middle Value Theorem, for some $\dot{\theta}$ between $\tilde{\theta}$ and $\theta^*$.
As usual, under the usual regularity conditions, since $\bigl(\tilde{\psi},\dot{\theta}\bigr)\rightarrow^P \bigl(\psi^*,\theta^*\bigr)$,
\begin{equation}\label{doublederlast}
\left.\frac{\partial}{\partial\theta}\right|_{\dot{\theta}}\frac{\partial}{\partial\theta_p}\mathbb{P}_n U_q\bigl(\tilde{\psi},\theta\bigr)
\rightarrow^P E\left.\frac{\partial}{\partial\theta}\right|_{\theta^*}\frac{\partial}{\partial\theta_p} U_q\bigl(\psi^*,\theta\bigr)=0,
\end{equation}
where the equality follows from the same reasoning as equation~(\ref{DRderivative}).
Combining equations~(\ref{doublederlast}) and~(\ref{noneforth})  implies that $n^{1/4}$ times the last term in equation~(\ref{doubleder}) converges in probability to zero. 

For the first term on the right hand side of equation~(\ref{doubleder}), because of the Middle Value Theorem,
\begin{eqnarray}\label{DRconswithnoneforth}
\lefteqn{\left.\frac{\partial}{\partial\theta_p}\right|_{\theta^*}\mathbb{P}_n U_q\bigl(\tilde{\psi},\theta\bigr)}\nonumber\\
&=&\left.\frac{\partial}{\partial\theta_p}\right|_{\theta^*}\mathbb{P}_n U_q(\psi^*,\theta)+\left(\left.\frac{\partial}{\partial\psi}\right|_{\dot{\psi}}\left.\frac{\partial}{\partial\theta_p}\right|_{\theta^*}\mathbb{P}_n U_q\bigl(\psi,\theta\bigr)\right)\bigl(\tilde{\psi}-\psi^*\bigr),
\end{eqnarray}
for some $\dot{\psi}$ between $\tilde{\psi}$ and $\psi$, possibly different in each row.
As usual, under the usual regularity conditions,
\begin{equation}\label{ddpsiddtheta}\left.\frac{\partial}{\partial\psi}\right|_{\dot{\psi}}\left.\frac{\partial}{\partial\theta_p}\right|_{\theta^*}\mathbb{P}_n U_q\bigl(\psi,\theta\bigr)\rightarrow^P E\left.\frac{\partial}{\partial\psi}\right|_{\psi^*}\left.\frac{\partial}{\partial\theta_p}\right|_{\theta^*} U_q\bigl(\psi,\theta\bigr).
\end{equation}
Moreover, from equation~(\ref{DRderivative}), the Central Limit Theorem implies that
\begin{equation}\label{cltonp}\sqrt{n}\left.\frac{\partial}{\partial\theta_p}\right|_{\theta^*}\mathbb{P}_n U_q(\psi^*,\theta)\rightarrow^{\cal D}{\cal N}\left(0,E\left(\left(\left.\frac{\partial}{\partial\theta_p}\right|_{\theta^*}U_q(\psi^*,\theta)\right)^{2}\right)\right).
\end{equation}
Combining equations~(\ref{DRconswithnoneforth}), (\ref{ddpsiddtheta}), (\ref{cltonp}), and (\ref{psinforth}) leads to
\begin{equation*}
n^{1/4}\left.\frac{\partial}{\partial\theta_p}\right|_{\theta^*}\mathbb{P}_n U_q\bigl(\tilde{\psi},\theta\bigr)\rightarrow^P0.
\end{equation*}
Combining with equations~(\ref{doubleder}), (\ref{doublederlast}), and (\ref{DRconswithnoneforth}), it follows that
\begin{equation}\label{almostthere}
n^{1/4}\left.\frac{\partial}{\partial\theta_p}\right|_{\tilde{\theta}}\mathbb{P}_n U_q\bigl(\tilde{\psi},\theta\bigr)\rightarrow^P 0.
\end{equation}

Combining equations~(\ref{der}), (\ref{almostthere}), and~(\ref{noneforth}), it follows that $\sqrt{n}$ times the last term in equation~(\ref{psiminpsistar}) converges in probability to zero. 

\section{Conclusion}

It follows that if equation~(\ref{noneforth}) and the usual regularity conditions from Section~\ref{regconditions} hold, for $\hat{\psi}$ of the form of Section~\ref{psisection} equation~(\ref{psihat}),
\begin{eqnarray}
\sqrt{n}\bigl(\hat{\psi}-\psi^*\bigr)
&=&\left(\left.\frac{\partial}{\partial\psi}\right|_{\tilde{\psi}}\mathbb{P}_n U\bigl(\psi,\tilde{\theta}\bigr)\right)^{-1} \sqrt{n}\mathbb{P}_n U\bigl(\psi^*,\theta^*\bigr)+o_P(1)\nonumber\\
&\rightarrow^{\cal D} &{\cal N}\left(0,V\left(\psi^*,\theta^*\right)\right),\label{result}
\end{eqnarray}
from the Central Limit Theorem, equation~(\ref{der}), and Slutsky's Theorem, with
\begin{equation}
V\left(\psi^*,\theta^*\right)=\left(E\left.\frac{\partial}{\partial\psi}\right|_{\psi^*}U\bigl(\psi,\theta^*\bigr)\right)^{-1}
E\left(U\bigl(\psi^*,\theta^*\bigr)^2\right)\left(E\left.\frac{\partial}{\partial\psi}\right|_{\psi^*}U\bigl(\psi,\theta^*\bigr)\right)^{-1\top}.
\end{equation}
That is, estimating $\theta$ leads to the same variance as plugging in the true but usually unknown $\theta$, and the sandwich estimator for the variance of $\hat{\psi}$ ignoring estimation of $\theta$ is consistent, all provided that $\hat{\psi}$ is doubly robust and $\theta$ is estimated at rate $n^{1/4}$ or faster.

\bibliographystyle{chicago}
\bibliography{ref}

\end{document}